\documentclass[11pt]{article}
\usepackage{graphicx}
\usepackage{color}
\usepackage{amsmath}
\usepackage{amssymb}
\usepackage{amscd}
\usepackage{framed}
\usepackage{bbm}

\usepackage{fancyhdr,a4wide}
\usepackage{amsthm}

\newcommand{\R}{\mathbb{R}}

\newcommand{\E}{\mathbb{E}}

\newtheorem{Theorem}{Theorem}[section]
\newtheorem{Lemma}[Theorem]{Lemma}

\newtheorem{Definition}[Theorem]{Definition}

\newtheorem{Remark}[Theorem]{Remark}

\numberwithin{equation}{section}

\def \proof {\noindent {\bf Proof.}\ \ }

\def \endproof
{{\mbox{}\nolinebreak\hfill\rule{2mm}{2mm}\par\medbreak}}

\newcommand{\wt}{\widetilde}
\newcommand{\wh}{\widehat}

\newcommand{\X}{\mathcal{X}}

\newcommand{\F}{{\mathcal F}}

\newcommand{\cX}{{\mathcal X}}

\begin{document}
\title{A remark on ``Robust machine learning by median-of-means"
%\thanks{
%G\'abor Lugosi was supported by
%the Spanish Ministry of Economy and Competitiveness,
%Grant MTM2015-67304-P and FEDER, EU. Shahar Mendelson was supported in part by the Israel Science Foundation.
%}
}
\author{
G\'abor Lugosi\thanks{ICREA and Pompeu Fabra University, gabor.lugosi@upf.edu}
\and
Shahar Mendelson \thanks{Department of Mathematics, Technion, I.I.T, and Mathematical Sciences Institute, The Australian National University, shahar@tx.technion.ac.il}}

\maketitle

\section{Introduction}

In this note we address regression function estimation, one of the most basic problems of statistical
learning. Suppose $\cX$ is some measurable space.
Given an independent, identically distributed sample
$D= (X_i,Y_i)_{i=1}^N$ of pairs of random variables with $X_i\in \cX$ and
$Y_i\in\R$, one wishes to select a function $\wh{f}: \X\to \R$ from a
given class $F$ with small risk
\[
\E \left( (\wh{f}(X)-Y)^2 | D \right)~.
\]
Tournament procedures were introduced in \cite{LugMen1} and attain the optimal accuracy/confidence tradeoff in general prediction problems (see \cite{LugMen1,LugMen2,Men}). Roughly put, the idea behind tournaments is to select $n \leq N$ wisely, split the given sample $(X_i,Y_i)_{i=1}^N$ to $n$ blocks, each of cardinality $m=N/n$, and compare the statistical performance of every pair of functions on each block. The function that exhibits a superior performance on the majority of the blocks is the winner of the match, and in a perfect world, the procedure returns a function that wins all of its matches. However, the world is far from perfect: the outcome of a match between functions that are too close is not reliable. To address this issue, tournaments require an additional component: a data-dependent way of verifying when two functions are too close, allowing one to decide if the result of a match should be trusted.

Although tournaments attain the optimal accuracy/confidence tradeoff, they are far from being computationally feasible: comparing every pair of functions in a large, possibly infinite class is impossible. A natural question is therefore to find a more reasonable procedure that exhibits the same optimal statistical behaviour as tournaments and at the same time has at least a fighting chance of being computationally friendly. The authors of \cite{LecLer} claim that the procedure they suggest has these features and those claims are the subject of this note.

Before we explore the claims from \cite{LecLer}, let us describe some technical facts that are at the heart of the results in \cite{LugMen1,LugMen2} and that will prove to be significant in what follows.

Given a class of functions $F$, let $f^*={\rm argmin}_{f \in F} \E(f(X)-Y)^2$. Fix an integer $n$ and let $(I_j)_{j=1}^n$ be the natural decomposition of $\{1,\ldots,N\}$ of $n$ blocks of cardinality $m=N/n$.  Set
$$
\mathbbm{Q}_{f,h}(j) = \frac{1}{m}\sum_{i \in I_j} (f-h)^2(X_i), \ \ \ \mathbbm{M}_{f,h}(j)= \frac{2}{m}\sum_{i \in I_j} (f-h)(X_i) (h(X_i)-Y_i)
$$
and note that
$$
\mathbbm{B}_{f,h}(j) = \frac{1}{m}\sum_{i \in I_j}(f(X_i)-Y_i)^2 - \frac{1}{m}\sum_{i \in I_j}(h(X_i)-Y_i)^2 = \mathbbm{Q}_{f,h}(j)+\mathbbm{M}_{f,h}(j)~.
$$

The method developed in \cite{LugMen1} was used there to show the
following fact: given a convex class $F$ that satisfies some minimal
conditions, for the right choice of $n$ and $r$ (the choice of $r$
depends on the geometry of the class $F$ and on the parameters
$\gamma_1$ and $\gamma_2$ appearing below), and for an absolute
constant $c_1$, we have that, with probability $1-2\exp(-c_1n)$,
\begin{framed}
\begin{description}
\item{(1)} for every $f \in F$ such that $\|f-f^*\|_{L_2} \geq r$, one has
\begin{equation} \label{eq:basic-1}
\mathbbm{B}_{f,f^*}(j) \geq \gamma_1 \|f-f^*\|_{L_2}^2
\end{equation}
for $0.99n$ of the blocks;
\item{(2)} for every $f \in F$ such that $\|f-f^*\|_{L_2} < r$, one has
\begin{equation} \label{eq:basic-2}
|\mathbbm{M}_{f,f^*}(j) - \E \mathbbm{M}_{f,f^*}(j)| \leq \gamma_2 r^2
\end{equation}
for $0.99n$ of the blocks.
\end{description}
Note that any fixed proportion of $n$ is possible, for the price of slightly modified constants.
\end{framed}
\eqref{eq:basic-1} and \eqref{eq:basic-2} are instrumental in our
analysis of the performance of the procedure from \cite{LecLer}.

\vskip0.4cm

The fact that \eqref{eq:basic-1} and \eqref{eq:basic-2} hold with probability at least $1-2\exp(-c_1n)$ is a reformulation of statements from \cite{LugMen1}: $(1)$ is Lemma 5.1 in \cite{LugMen1} and $(2)$ is Lemma 5.5 in \cite{LugMen1}.

\begin{Remark}
Note that both facts are not enough to run a tournament procedure: although the identity of the winner of each match involving $f^*$ is clear, it does not tell us whether the outcome of the match should be trusted. However, using the right notion of a distance oracle, these estimates suffice to ensure that with probability at least $1-2\exp(-c n)$ the tournament procedure produces $\wh{f}$ for which
$$
\E \left( (\wh{f}(X)-Y)^2 | (X_i,Y_i)_{i=1}^N \right)\leq \E (f^*(X)-Y)^2 + c^\prime r^2.
$$
\end{Remark}

As mentioned previously, the motivation for this note is \cite{LecLer}, where the authors suggest the following alternative to the tournament procedure: given a convex class $F$ and a sample $(X_i,Y_i)_{i=1}^N$, select $\wt{f}$ to be the minimizer in $F$ of the functional
$$
\phi(f) = \max_{g \in F} {\rm Med}( \mathbbm{B}_{f,g}),
$$
where ${\rm Med}( \mathbbm{B}_{f,g})$ is a median of the vector $(\mathbbm{B}_{f,g}(j))_{j=1}^n$.
\vskip0.5cm
The claim in \cite{LecLer} is that for a well chosen number of blocks $n$ (the same as in \cite{LugMen1,LugMen2}) $\wt{f}$ performs as well as the tournament procedure---both in the standard prediction framework of \cite{LugMen1} and in the regularization framework of \cite{LugMen2}; that it is robust to outliers and that it is computationally feasible.

To prove their claims, the authors make the same assumptions as in \cite{LugMen1,LugMen2}, and essentially re-prove the technical machinery developed in \cite{LugMen1}. This machinery is then used to analyze the performance of $\wt{f}$.

\vskip0.3cm
We show below that the theoretical contribution of \cite{LecLer} is somewhat overstated: the two main claims on the statistical performance of $\wt{f}$ are in fact almost obvious outcomes of \eqref{eq:basic-1} and \eqref{eq:basic-2}, and thus the analysis of $\wt{f}$ is a direct and immediate corollary of Lemma 5.1 and Lemma 5.5 from \cite{LugMen1}. Moreover, we show that the robustness to outliers exhibited by $\wt{f}$ is equally simple (this is one of the features of quantile-based procedures and therefore should not come as a major surprise).
\vskip0.3cm
Finally, \cite{LecLer} explores the issue of computational feasibility of $\wt{f}$. The claim is that one may use coordinate descent to find an approximate solution to the minimization problem that defines $\wt{f}$. We believe that at this point it is yet to be determined whether coordinate descent truly finds an approximate solution in polynomial time with the same theoretical properties as $\wt{f}$. Specifically, we question whether the optimal tradeoff between accuracy and confidence---the main novelty in tournaments---can be ensured by such an approximate solution.

\section{The analysis of $\wt{f}$}
Assume that one is given a sample $D=(X_i,Y_i)_{i=1}^N$ for which \eqref{eq:basic-1} and \eqref{eq:basic-2} hold.
The analysis of the procedure $\wt{f}$ is straightforward: observe that as a minimizer, $\phi(\wt{f}) \leq \phi(f^*)$, and moreover, $\phi(f^*) = \max_{g \in F} {\rm Med}(-\mathbbm{B}_{g,f^*})$. Thanks to \eqref{eq:basic-1} and \eqref{eq:basic-2} one may provide an upper bound on $\phi(f^*)$ and then use it to pin-point the location of $\wt{f}$ in $F$.

\begin{Theorem} \label{thm:main1}
Assume that $\gamma_1>\gamma_2$. Then for the given sample $D$ one has
$$
\E \left( \left(\wt{f}(X)-Y\right)^2|D \right) \leq \E (f^*(X)-Y)^2 + (1+2\gamma_2)r^2.
$$
\end{Theorem}

\proof
Let $g \in F$ satisfy that $\|g-f^*\|_{L_2} \geq r$. By \eqref{eq:basic-1}, $-\mathbbm{B}_{g,f^*}(j) \leq -\gamma_1 r^2$ on $0.99n$ of the blocks. On the other hand, if $\|g-f^*\|_{L_2} \leq r$, we have that
$$
\mathbbm{B}_{g,f^*}(j) = \mathbbm{Q}_{g,f^*}(j) + \mathbbm{M}_{g,f^*}(j) \geq \mathbbm{M}_{g,f^*}(j).
$$
Since $F$ is convex one has that $\E \mathbbm{M}_{g,f^*}(j) \geq 0$. Indeed, this follows from the characterization of the nearest point map onto a convex set in a Hilbert space. Hence,
by \eqref{eq:basic-2},
\begin{equation} \label{eq:multi1}
\mathbbm{M}_{g,f^*}(j) \geq \E \mathbbm{M}_{g,f^*}(j) - \gamma_2 r^2 \geq - \gamma_2 r^2
\end{equation}
on $0.99n$ of the blocks.

Combining these two observations,
$$
\phi(f^*) = \max_{g \in F} {\rm Med}(-\mathbbm{B}_{g,f^*}) \leq  \gamma_2 r^2.
$$
Now, as the minimizer, $\phi(\wt{f}) = \max_{g \in F} {\rm
  Med}(\mathbbm{B}_{\wt{f},g}) \leq \gamma_2 r^2$. In particular,
${\rm Med}(\mathbbm{B}_{\wt{f},f^*}) \leq \gamma_2 r^2$, which forces
$\wt{f}$ to be in a rather specific part of $F$. To identify the
location of $\wt{f}$, first observe that $\|\wt{f}-f^*\|_{L_2} \leq
r$. Otherwise, by \eqref{eq:basic-1}, $\mathbbm{B}_{\wt{f},f^*} \geq
\gamma_1 r^2$ on $0.99n$ of the  blocks, which is impossible if
$\gamma_1 > \gamma_2$. Finally, given that $\|\wt{f}-f^*\|_{L_2} \leq
r$ then also $\E \mathbbm{M}_{\wt{f},f^*} \leq 2\gamma_2 r^2$, since
otherwise, by \eqref{eq:multi1},
we would have
$$
\mathbbm{B}_{\wt{f},f^*}(j) \geq \mathbbm{M}_{\wt{f},f^*}(j) \geq \E \mathbbm{M}_{\wt{f},f^*}(j) - \gamma_2 r^2 > \gamma_2 r^2
$$
on $0.99n$ of the blocks, which is also impossible.

Therefore we have that
$$
\E \left( (\wt{f}(X)-Y)^2 | D \right) - \E (f^*(X)-Y)^2 = \|\wt{f}-f^*\|_{L_2}^2 + \E \mathbbm{M}_{\wt{f},f^*} \leq (1+2\gamma_2) r^2,
$$
as claimed.
\endproof

\begin{Remark}
The second claim from \cite{LecLer}, that just like tournaments, $\wt{f}$ is robust to malicious corruption of the given data, is now completely clear: if \eqref{eq:basic-1} and \eqref{eq:basic-2} are to be believed, it means that $98\%$  of the values of both $\mathbbm{B}_{f,f^*}(j)$ and $\mathbbm{M}_{f,f^*}(j)$ are in the range we want. Even if another $40\%$ of those values are changed maliciously, the median of $(\mathbbm{B}_{f,f^*}(j))_{i=1}^n$ would still be in the right range, and the proof of Theorem \ref{thm:main1} would still hold. As noted previously, this robustness to malicious changes is one of main features of quantiles.

It should be noted that the assumption made in \cite{LecLer} is that the number of ``corrupted samples'' is smaller than a small proportion of the number of blocks, as one would expect from the proof of Theorem \ref{thm:main1}.
\end{Remark}

\subsection{Regularization}
Next, let us turn to the question of regularization, where the analysis of $\wt{f}$ is equally simple and again, requires only the right versions of Lemma 5.1 and Lemma 5.5 from \cite{LugMen1}.

In the regularized tournament introduced in \cite{LugMen2}, the match between $f$ and $h$ is determined by the $n$ values
\begin{equation} \label{eq:regularization1}
\mathbbm{B}_{f,h}(j) + \lambda (\Psi(f) - \Psi(h) ),
\end{equation}
where $\Psi$ denotes the regularization function and $\lambda$ is the regularization parameter. The alternative version of the regularized tournament suggested in \cite{LecLer} was to select $\wt{f}_\lambda$, set to be the minimizer of the functional
\begin{equation} \label{eq:regularization2}
\phi_\lambda(f) = \max_{g \in F} {\rm Med} (\mathbbm{B}_{f,g}) + \lambda (\Psi(f)-\Psi(g)).
\end{equation}
Before we explain how the performance of $\wt{f}_\lambda$ can be determined using \eqref{eq:basic-1} and \eqref{eq:basic-2} and in order to put our presentation in some context, let us briefly describe some ideas that are needed in the study of regularized procedures. For a more detailed description we refer the reader to \cite{LecMen}.

The whole point in regularization is dealing with a problem that involves a convex class $F$ that is simply `too big'. To address the size of $F$, one assigns a `price-tag' to each function in the class according to some prior belief, giving less favourable functions a higher price. That price is captured by the regularization function $\Psi$, which, for the sake of simplicity is assumed to be a norm on a linear space $E$ that contains $F$.

As always, one would like to use the given random data $(X_i,Y_i)_{i=1}^N$ to `distinguish' between functions $f \in F$ and $f^*$, allowing one to rule-out functions whose statistical performance is inferior to that of $f^*$. A good procedure is designed to select only functions that cannot be `distinguished' from $f^*$, and the regularized procedures we consider are no exception.

The method of analysis we use here is based on ideas from \cite{LecMen}: given the class $F$, we first identify two radii, $\rho$ and $r$; consider ${\cal B}_{f^*}(\rho)$---the $\Psi$-ball centred at $f^*$ and of radius $\rho$, and set $F_\rho = F \cap {\cal B}_{f^*}(\rho)$. Now one proceeds with the following three steps:
\begin{description}
\item{(1)} Since the set $F_\rho$ is much smaller than $F$, standard
  methods of `distinguishing' between $f^*$ and functions in $F_\rho$
  is possible. Moreover, the difference between $\Psi(f)$ and
  $\Psi(f^*)$ is not that big---at most $\rho$. Hence, one would like to show that for functions in $F_\rho$ that satisfy  $\|f-f^*\|_{L_2} \geq r$, the empirical component of the excess functional is `positive enough' to overcome the contribution of the regularization terms, which is no worse than $-\lambda \rho$.
\item{(2)} When considering functions in $F$ that satisfy $\Psi(f-f^*)=\rho$ and for which $\|f-f^*\|_{L_2} \leq r$ it is futile to expect that the empirical component of the excess functional is of any use. This is precisely when functions are `too close' and cannot be distinguished using empirical data. Therefore, all the `hard work' required to distinguish such functions from $f^*$ has to be based on properties of the regularization function $\Psi$. We describe how to obtain the wanted control in what follows.
\item{(3)} Combining $(1)$ and $(2)$ we can distinguish between $f^*$ and any function that satisfies $\Psi(f-f^*)=\rho$, or between $f^*$ and functions in $F_\rho$ that satisfy $\|f-f^*\|_{L_2} \geq r$. Moreover, by a homogeneity argument, the estimate in the $\Psi$-sphere transfers for free to functions $f \in F$ that satisfy $\Psi(f-f^*)>\rho$, allowing us to distinguish between those and $f^*$.
\end{description}
Using the combination of the three components one may show that a procedure selects $h$ that satisfies both $\Psi(h-f^*) \leq \rho$ and $\|h-f^*\|_{L_2} \leq r$; showing that in addition, the excess risk of $h$ is smaller than $cr^2$ requires an additional argument, and we return to it later.

\vskip0.5cm

The key component in regularization happens to be $(2)$: identifying features of the regularization function $\Psi$ that yield sufficient control when the empirical component of the functional fails. By now it is well understood that this property has to do with the smoothness of $\Psi$, as we briefly explain next.

Let $B_{\Psi^*}$ and $S_{\Psi^*}$ denote the unit ball and unit sphere in the dual space to $(E,\Psi)$, respectively. Therefore, $B_{\Psi^*}$ consists of all the linear functionals $z \in E^*$ for which $\sup_{\{x \in E : \Psi(x)=1\}} |z(x)| \leq 1$. A linear functional $z^* \in S_{\Psi^*}$ is a norming functional for $f \in E$ if $z^*(f)=\Psi(f)$.

\begin{Definition} \label{def:Delta}
Let $\Gamma_{f}(\rho) \subset S_{\Psi^*}$ be the collection of functionals that are norming for some $v \in {\cal B}_f(\rho/20)$. Set
\begin{equation*}
\Delta_F(\rho,r) = \inf_{f \in F} \inf_h \sup_{z \in \Gamma_{f}(\rho)} z(h-f)~,
\end{equation*}
where the inner infimum is taken in the set
\begin{equation} \label{eq:condition-on-set}
\{ h \in F: \Psi(h-f)=\rho \ {\rm and} \ \|h-f\|_{L_2} \leq r\}~.
\end{equation}
\end{Definition}

Several examples of regularization functions $\Psi$ and the resulting estimates on $\Delta_F(\rho,r)$ can be found in \cite{LecMen,LugMen2}. Among the examples are standard sparsity-driven procedures like {\sc lasso} and {\sc slope}.

\vskip0.4cm
For our purposes, the crucial observation incorporates a lower bound on $\Delta_F(\rho,r)$ and a wise choice of the regularization parameter $\lambda$. This observation is not new as well: it is a version of Lemma 4.5 from \cite{LugMen2} and although its formulation is not identical to that lemma, its proof is --- line for line. We present the proof in an the appendix for the sake of completeness.

\begin{Lemma} \label{lemma:basic-combining-loss-and-reg}
Fix $\gamma_1$ and $\gamma_2$ and some  block $I_j$. Let $\rho$ and $r$ such that $\Delta_F(\rho,r) \geq 4\rho/5$ and set $\lambda$ to satisfy
\begin{equation} \label{eq:lambda}
3\gamma_2 \cdot \frac{r^2}{\rho} \leq \lambda \leq \frac{\gamma_1}{2} \cdot \frac{r^2}{\rho}.
\end{equation}
Assume that for $h \in F_\rho$ such that $\|h-f^*\|_{L_2} \geq r$, one has
\begin{equation} \label{eq:cond3}
\mathbbm{B}_{h,f^*}(j) \geq \gamma_1 \|h-f^*\|_{L_2}^2;
\end{equation}
and if $\|h-f^*\|_{L_2} \leq r$ then
\begin{equation} \label{eq:cond4}
|\mathbbm{M}_{h,f^*}(j)-\E \mathbbm{M}_{h,f^*}(j)| \leq \gamma_2 r^2 .
\end{equation}
Under these conditions the following holds:
\begin{description}
\item{(1)} If $h \in F_\rho$ and $\|h-f^*\|_{L_2} \geq r$ then
\begin{equation} \label{eq:reg-control-1}
\mathbbm{B}_{h,f^*}(j) + \lambda(\Psi(h)-\Psi(f^*)) \geq \frac{\gamma_1}{2}\|h-f^*\|_{L_2}^2;
\end{equation}
\item{(2)} If $\Psi(h-f^*)=\rho$ and $\|h-f^*\|_{L_2} < r$ then
\begin{equation} \label{eq:reg-control-2}
\mathbbm{B}_{h,f^*}(j) + \lambda(\Psi(h)-\Psi(f^*)) \geq  \frac{\gamma_2}{2} r^2;
\end{equation}
\item{(3)} Let $f,h \in F$ satisfy that $\Psi(h-f^*)=\rho$ and $f=f^*+\alpha(h-f^*)$ for some $\alpha>1$. If $\|h-f^*\|_{L_2} \geq r$ then
$$
\mathbbm{B}_{f,f^*}(j) + \lambda(\Psi(f)-\Psi(f^*)) \geq  \alpha \cdot \frac{\gamma_1}{2}\|h-f^*\|_{L_2}^2,
$$
and if $\|h-f^*\|_{L_2} < r$ then
$$
\mathbbm{B}_{f,f^*}(j) + \lambda(\Psi(f)-\Psi(f^*)) \geq  \alpha \cdot \gamma_2 r^2
$$
\end{description}
\end{Lemma}

Conditions \eqref{eq:cond3} and \eqref{eq:cond4} should not come as a surprise: they are simply \eqref{eq:basic-1} and \eqref{eq:basic-2} for the class $F_\rho$ and for one block. Hence, by Lemma 5.1 and Lemma 5.5 from \cite{LugMen1}, under minimal assumptions, with probability at least $1-2\exp(-c_1n)$ for every $h \in F_\rho$, each condition holds for $0.99n$ of the blocks.

\begin{Theorem} \label{thm:main2}
Let $r,\rho$ and $\lambda$ be as in Lemma \ref{lemma:basic-combining-loss-and-reg} and assume that for a sample $D=(X_i,Y_i)_{i=1}^N$, \eqref{eq:basic-1} and \eqref{eq:basic-2} hold in the class $F_\rho$. Then,
$$
\Psi(\wt{f}_\lambda-f^*) \leq c_1\rho, \ \ \ \ \|\wt{f}_\lambda-f^*\|_{L_2} \leq c_2r
$$
and
$$
\E \left( (\wt{f}_\lambda(X)-Y)^2 | D \right) - \E (f^*(X)-Y)^2 \leq c_3r^2
$$
where $c_1,c_2$ and $c_3$ depend only on $\gamma_1$ and $\gamma_2$.
\end{Theorem}

\proof Recall the properties $(1)$, $(2)$ and $(3)$ from Lemma \ref{lemma:basic-combining-loss-and-reg}, which we first use to obtain an upper estimate on
$$
\phi_\lambda(f^*) = \max_{g \in F} \left[ - \left( {\rm Med} (\mathbbm{B}_{g,f^*}) + \lambda (\Psi(g)-\Psi(f^*))\right) \right]
$$
by exploring all the possible options of $g \in F$.

If $g \in F_\rho$ and $\|g-f^*\|_{L_2} \geq r$ then by $(1)$, for $0.99n$ of the blocks,
$$
\mathbbm{B}_{g,f^*}(j) + \lambda(\Psi(g)-\Psi(f^*)) \geq \frac{\gamma_1}{2}\|g-f^*\|_{L_2}^2;
$$
and if $\Psi(g-f^*)=\rho$ and $\|g-f^*\|_{L_2} < r$ then by $(2)$, on $0.99n$ of the blocks,
$$
\mathbbm{B}_{g,f^*}(j) + \lambda(\Psi(g)-\Psi(f^*)) \geq \frac{\gamma_2}{2} r^2.
$$
Thus, together with the ``super-linear" growth of $\mathbbm{B}_{h,f^*}(j) + \lambda(\Psi(h)-\Psi(f^*))$ from $(3)$, we have that if $\Psi(g-f^*) \geq \rho$ then
$$
- \left( {\rm Med} (\mathbbm{B}_{g,f^*}) + \lambda (\Psi(g)-\Psi(f^*))\right) < 0.
$$
All that remains is to study the case $\Psi(g-f^*) < \rho$ and $\|g-f^*\|_{L_2} \leq r$. Note that in that range, for $0.99n$ of the blocks, $|\mathbbm{M}_{g,f^*}(j)-\E \mathbbm{M}_{g,f^*}(j)| \leq \gamma_2 r^2 $ and by the convexity of $F$, $\E \mathbbm{M}_{g,f^*}(j)>0$. Therefore,
$$
\mathbbm{B}_{g,f^*}(j) \geq \mathbbm{M}_{g,f^*}(j) \geq -\gamma_2 r^2 + \E \mathbbm{M}_{g,f^*}(j) \geq -\gamma_2 r^2,
$$
and with our choice of $\lambda$,
$$
\mathbbm{B}_{g,f^*}(j) + \lambda (\Psi(g)-\Psi(f^*)) \geq \gamma_2 r^2 - \lambda \rho  \geq -\left(\gamma_2 + \frac{\gamma_1}{2}\right)r^2.
$$
Hence, we have that $\phi_\lambda(f^*) \leq \left(\gamma_2 + \frac{\gamma_1}{2}\right)r^2$, and just as we did previously, we use this information to pin-point the location of $\wt{f}_\lambda$.

By considering the choice $g=f^*$, we have
$$
\phi_\lambda(\wt{f}_\lambda) \geq {\rm Med}(\mathbbm{B}_{\wt{f}_\lambda,f^*}) + \lambda(\Psi(\wt{f}_\lambda)-\Psi(f^*)),
$$
and let us rule out possible locations of $\wt{f}_\lambda$. To that end, let $\alpha>1$ to be specified later and set $f \in F$ such that $\Psi(f-f^*) = \alpha \rho$, i.e., $f=f^*+\alpha(h-f^*)$, where $h \in F$ and $\Psi(h-f^*)=\rho$. By the super-linear growth in $(3)$, combined with $(1)$ and $(2)$, it follows that for $0.98n$ of the blocks,
\begin{equation*}
\mathbbm{B}_{f,f^*}(j) + \lambda(\Psi(f)-\Psi(f^*)) \geq \frac{\alpha}{2} \min\left\{\gamma_1,\gamma_2\right\} r^2 > \left(\gamma_2 + \frac{\gamma_1}{2}\right)r^2
\end{equation*}
for the right choice of a large enough $\alpha$ that depends only on $\gamma_1$ and $\gamma_2$. Repeating this argument for larger values $\alpha^\prime > \alpha $ rules out the possibility that $\Psi(\wt{f}_\lambda-f^*) \geq \alpha \rho$, implying that $\Psi(\wt{f}_\lambda-f^*) < \alpha \rho$.

Given that $\Psi(\wt{f}_\lambda-f^*) \leq \alpha \rho$, let us estimate $\|\wt{f}-f^*\|_{L_2}$. If $\|\wt{f}_\lambda-f^*\|_{L_2} \geq \alpha r$, set
\begin{equation} \label{eq:h}
h=f^*+ \frac{1}{\alpha}(\wt{f}_\lambda-f^*)
\end{equation}
and observe that $\Psi(h-f^*) \leq \rho$ and $\|h-f^*\|_{L_2} \geq r$. Hence, by $(1)$, on $0.99n$ of the blocks,
$$
\mathbbm{B}_{h,f^*}(j) + \lambda(\Psi(h)-\Psi(f^*)) \geq \frac{\gamma_1}{2}\|h-f^*\|_{L_2}^2,
$$
and by $(3)$ we have for the right choice of $\alpha$ that on the same blocks
\begin{equation*}
\mathbbm{B}_{\wt{f}_\lambda,f^*}(j) + \lambda(\Psi(\wt{f}_\lambda)-\Psi(f^*)) \geq \alpha \frac{\gamma_1}{2}\|h-f^*\|_{L_2}^2
\geq \alpha \frac{\gamma_1}{2} r^2 > \left(\gamma_2 + \frac{\gamma_1}{2}\right)r^2,
\end{equation*}
which is impossible. Hence, $\Psi(\wt{f}_\lambda-f^*) \leq \alpha \rho$ and $\|\wt{f}_\lambda-f^*\|_{L_2} \leq \alpha r$. Moreover, $h$ defined in \eqref{eq:h} satisfies that $\|h-f^*\|_{L_2} \leq r$.

Finally, let us show that
\begin{equation} \label{eq:in-reg-M}
\E \mathbbm{M}_{\wt{f}_\lambda,f^*}(j) \leq 2(\alpha+1)\left(\gamma_2 + \frac{\gamma_1}{2}\right)r^2.
\end{equation}
Indeed, assume that the reverse inequality holds. Since $\|h-f^*\|_{L_2} \leq r$ then by \eqref{eq:cond4}, on $0.99n$ of the blocks, $|\mathbbm{M}_{h,f^*}(j)-\E \mathbbm{M}_{h,f^*}(j)| \leq \gamma_2 r^2$; and, since $\alpha \mathbbm{M}_{h,f^*}(j) = \mathbbm{M}_{\wt{f}_\lambda,f^*}(j)$ it follows that on the same blocks,
$$
\mathbbm{M}_{\wt{f}_\lambda,f^*}(j) \geq \E \mathbbm{M}_{\wt{f}_\lambda,f^*}(j) - \alpha \gamma_2 r^2.
$$
Thus, by our choice of $\lambda$,
\begin{align*}
& \mathbbm{B}_{\wt{f}_\lambda,f^*}(j) + \lambda (\Psi(\wt{f}_\lambda)-\Psi(f^*)) \geq
\E \mathbbm{M}_{\wt{f}_\lambda,f^*}(j) - \alpha \gamma_2 r^2 - \lambda \cdot \alpha \rho
\\
\geq & \E \mathbbm{M}_{\wt{f}_\lambda,f^*}(j) - \alpha \left(\gamma_2+\frac{\gamma_1}{2}\right) r^2 >\left(\gamma_2 + \frac{\gamma_1}{2}\right)r^2,
\end{align*}
which is impossible --- confirming \eqref{eq:in-reg-M}.

To conclude, we have shown that for the sample $D=(X_i,Y_i)_{i=1}^N$,

$$
\E \left( (\wt{f}_\lambda(X)-Y)^2 | D \right) - \E (f^*(X)-Y)^2 =
\|\wt{f}_\lambda-f^*\|_{L_2}^2 + \E \mathbbm{M}_{h,f^*} \leq c(\gamma_1,\gamma_2)r^2,
$$
completing the proof.
\endproof

\bibliographystyle{plain}
\bibliography{tournament3}

\newpage
\appendix

\section{Proof of Lemma \ref{lemma:basic-combining-loss-and-reg}}
Let us begin by examining
$$
(*) = \mathbbm{B}_{h,f^*}(j)+\lambda(\Psi(h)-\Psi(f^*))
$$
in the set $\{f \in \F : \Psi(f-f^*)=\rho\}$, where one should consider two cases. First, if $\|f-f^*\|_{L_2} \geq r$ and since $\mathbbm{B}_{h,f^*}(j) \geq \gamma_1 \|h-f^*\|_{L_2}^2$, then by the triangle inequality for $\Psi$, \begin{equation} \label{eq:large-L-2}
(*) \geq \gamma_1 \|h-f^*\|_{L_2}^2  - \lambda \Psi(f-f^*) \geq \gamma_1 \|h-f^*\|_{L_2}^2  - \lambda \rho  \geq \frac{\gamma_1}{2} \|h-f^*\|_{L_2}^2
\end{equation}
provided that
\begin{equation} \label{eq:lambda-upper}
\lambda \leq \frac{\gamma_1}{2} \cdot \frac{r^2}{\rho}~.
\end{equation}
If, on the other hand, $\|f-f^*\|_{L_2} \leq r$, then recalling that $\E \mathbbm{M}_{h,f^*}(j) \geq 0$ we have $\mathbbm{B}_{h,f^*}(j) \geq \mathbbm{M}_{h,f^*}(j) \geq -\gamma_2 r^2$; therefore, $(*) \geq -\gamma_2 r^2 + \lambda(\Psi(f)-\Psi(f^*))$.

Fix $v \in {\cal B}_{f^*}(\rho/20)$ and write $f^*=u+v$; thus $\Psi(u) \leq \rho/20$.  Set $z$ to be a linear functional that is norming for $v$ and observe that for any $h \in E$,
\begin{eqnarray} \label{eq:small-L-2a}
\Psi(h)-\Psi(f^*) & \geq & \Psi(h) - \Psi(v) - \Psi(u) \geq z(h-v)-\Psi(u) \geq z(h-f^*)-2\Psi(u) \nonumber
\\
& \geq & z(h-f^*) - \frac{\rho}{10}~.
\end{eqnarray}
Hence, if $f^* \in F$, $\Psi(h-f^*) =\rho$ and $\|h-f^*\|_{L_2} \leq r$,
then optimizing the choices of $v$ and of $z$, $z(h-f^*) \geq \Delta_F(\rho,r)$; thus
\begin{equation} \label{eq:small-L-2b}
\Psi(h)-\Psi(f^*) \geq \Delta_F(\rho,r) -\frac{\rho}{10} \geq \frac{7}{10} \rho~.
\end{equation}
And, if
\begin{equation} \label{eq:lambda-lower}
\lambda \geq 3\gamma_2 \cdot \frac{r^2}{\rho}~,
\end{equation}
we have that
$$
(*) \geq \frac{\gamma_2}{2} r^2.
$$

Next, if $h \in F_\rho$ and $\|h-f^*\|_{L_2} \geq r$, then
$$
\mathbbm{B}_{h,f^*}(j)+\lambda(\Psi(f)-\Psi(f^*)) \geq \frac{\gamma_1}{2} \|h-f^*\|_{L_2}^2;
$$
indeed, this follows from \eqref{eq:large-L-2}.

Finally, let us prove the super-linearity property when $\Psi(f-f^*) > \rho$. Set $\theta \in (0,1)$ and let $h \in F$ satisfy that
$$
\Psi(h-f^*) = \rho \ \ {\rm and} \ \ \theta (f-f^*) = h-f^*.
$$
If $\|h-f^*\|_{L_2} \geq r$, then by the triangle inequality for $\Psi$ followed by \eqref{eq:large-L-2},
\begin{eqnarray*}
(*) & \geq & \frac{1}{\theta^2} \mathbbm{Q}_{h,f^*}(j) + \frac{1}{\theta}\left(\mathbbm{M}_{h,f^*}(j) - \lambda \Psi(h-f^*) \right)
\\
& \geq & \frac{1}{\theta} \left(\mathbbm{B}_{h,f^*}(j) - \lambda \Psi(h-f^*) \right) \geq \frac{1}{\theta} \cdot \frac{\gamma_1}{2} \|h-f^*\|_{L_2}^2.
\end{eqnarray*}
If, on the other hand, $\|h-f^*\|_{L_2} \leq r$, then using the argument from \eqref{eq:small-L-2a}, \eqref{eq:small-L-2b} and the choice of $\lambda$ from \eqref{eq:lambda-lower}, we have
\begin{eqnarray*}
(*)
& \geq & \frac{1}{\theta} \left(\mathbbm{M}_{h,f^*}(j) +
  \lambda\left(z(h-f^*)-2\theta \Psi(u)\right)\right)
\\
& \geq & \frac{1}{\theta} \left(\mathbbm{M}_{h,f^*}(j) + \lambda\left(z(h-f^*)-2\Psi(u)\right)\right) \geq \frac{1}{\theta} \cdot \frac{\gamma_2}{2} r^2,
\end{eqnarray*}
as claimed.
\endproof

\end{document}